\documentclass{article}%
\usepackage{amsfonts}
\usepackage{amssymb}
\usepackage{amsmath}
\usepackage{graphicx}%
\providecommand{\U}[1]{\protect\rule{.1in}{.1in}}

\begin{document}

\title{On Adjoint and Brain Functors}
\author{David Ellerman\\Philosophy Department\\University of California at Riverside}
\date{}
\maketitle

\begin{abstract}
There is some consensus among orthodox category theorists that the concept of
adjoint functors is the most important concept contributed to mathematics by
category theory. We give a heterodox treatment of adjoints using
heteromorphisms (object-to-object morphisms between objects of different
categories) that parses an adjunction into two separate parts (left and right
representations of heteromorphisms). Then these separate parts can be
recombined in a new way to define a cognate concept, the brain functor, to
abstractly model the functions of perception and action of a brain. The
treatment uses relatively simple category theory and is focused on the
interpretation and application of the mathematical concepts. The Mathematical
Appendix is of general interest to category theorists as it is a defense of
the use of heteromorphisms as a natural and necessary part of category theory.

\end{abstract}
\tableofcontents

\section{Category theory in the life and cognitive sciences}

There is already a considerable but widely varying literature on the
application of category theory to the life and cognitive sciences--such as the
work of Robert Rosen (\cite{rosen:1958cat-theory}, \cite{rosen:anticipsys2ed})
and his followers\footnote{See \cite{zafiris:rosen},
\cite{louie:catsystemtheory}, and \cite{louie-poli:hspread} and their
references.} as well as Andr\'{e}e Ehresmann and Jean-Paul Vanbremeersch
\cite{ehresmann-v:mes} and their commentators.\footnote{See
\cite{kainen:mathbio} for Kainen's comments on the Ehresmann-Vanbremeersch
approach, Kainen's own approach, and a broad bibliography of relevant papers.}

The approach taken here is based on a specific use of the characteristic
concepts of category theory, namely universal mapping properties. One such
approach in the literature is that of Fran\c{c}ois Magnan and Gonzalo Reyes
which emphasizes that "Category theory provides means to circumscribe and
study what is universal in mathematics and other scientific disciplines."
\cite[p. 57]{mag-reyes:cog}. Their intended field of application is cognitive science.

\begin{quotation}
We may even suggest that universals of the mind may be expressed by means of
universal properties in the theory of categories and much of the work done up
to now in this area seems to bear out this suggestion....

By discussing the process of counting in some detail, we give evidence that
this universal ability of the human mind may be conveniently conceptualized in
terms of this theory of universals which is category theory. \cite[p.
59]{mag-reyes:cog}
\end{quotation}

Another current approach that emphasizes universal mapping properties
("universal constructions") is that of S. Phillips, W. H. Wilson, and G. S.
Halford (\cite{halford-wilson:cogdev}, \cite{phillips-wilson:uniconst},
\cite{phillips:analogy}).

In addition to the focus on universals, the approach here is distinctive in
the use of heteromorphisms--which are object-to-object morphisms between
objects if different categories--in contrast to the usual homomorphisms or
homs between objects in the same category. By explicitly adding
heteromorphisms to the usual homs-only presentation of category theory, this
approach can directly represent interactions between the objects of different
categories (intuitively, between an organism and the environment). But it is
still early days, and many approaches need to be tried to find out "where
theory lives."

\section{The ubiquity and importance of adjoints}

Before developing the concept of a brain functor, we need to consider the
related concept of a pair of adjoint functors, an adjunction. The developers
of category theory, Saunders Mac\thinspace Lane and Samuel Eilenberg, famously
said that categories were defined in order to define functors, and functors
were defined in order to define natural transformations
\cite{eilenberg-macl:gentheory}. A few years later, the concept of universal
constructions or universal mapping properties was isolated
(\cite{maclane:groups-duality} and \cite{samuel:universal}). Adjoints were
defined a decade later by Daniel Kan \cite{kan:adjoints} and the realization
of their ubiquity ("Adjoint functors arise everywhere" \cite[p. v]%
{maclane:cwm}) and their foundational importance has steadily increased over
time (Lawvere \cite{lawvere:adjointness} and Lambek \cite{lambek:heraclitus}).
Now it would perhaps not be too much of an exaggeration to see categories,
functors, and natural transformations as the prelude to defining adjoint
functors. As Steven Awodey put it:

\begin{quotation}
\noindent The notion of adjoint functor applies everything that we have
learned up to now to unify and subsume all the different universal mapping
properties that we have encountered, from free groups to limits to
exponentials. But more importantly, it also captures an important mathematical
phenomenon that is invisible without the lens of category theory. Indeed, I
will make the admittedly provocative claim that adjointness is a concept of
fundamental logical and mathematical importance that is not captured elsewhere
in mathematics. \cite[p. 179]{awodey:cat}
\end{quotation}

Other category theorists have given similar testimonials.

\begin{quotation}
\noindent To some, including this writer, adjunction is the most important
concept in category theory. \cite[p. 6]{wood:adjunctions}

\noindent The isolation and explication of the notion of adjointness is
perhaps the most profound contribution that category theory has made to the
history of general mathematical ideas." \cite[p. 438]{goldblatt:topoi}

\noindent Nowadays, every user of category theory agrees that [adjunction] is
the concept which justifies the fundamental position of the subject in
mathematics. \cite[p. 367]{taylor:foundations}
\end{quotation}

\section{Adjoints and universals}

How do the ubiquitous and important adjoint functors relate to the universal
constructions? Mac Lane and Birkhoff succinctly state the idea of the
universals of category theory and note that adjunctions can be analyzed in
terms of those universals.

\begin{quotation}
\noindent The construction of a new algebraic object will often solve a
specific problem in a universal way, in the sense that every other solution of
the given problem is obtained from this one by a unique homomorphism. The
basic idea of an adjoint functor arises from the analysis of such universals.
\cite[p. v]{maclane:algebra}
\end{quotation}

\noindent We can use some old language from Plato's theory of universals to
describe those universals of category theory (Ellerman
\cite{ellerman:erkenntnis}) that solve a problem in a universal or
paradigmatic way so that "every other solution of the given problem is
obtained from this one" in a unique way.

In Plato's Theory of Ideas or Forms ($\varepsilon\iota\delta\eta$), a property
$F$ has an entity associated with it, the universal $u_{F}$, which uniquely
represents the property. An object $x$ has the property $F$, i.e., $F(x)$, if
and only if (iff) the object $x$ \textit{participates} in the universal
$u_{F}$. Let $\mu$ (from $\mu\varepsilon\theta\varepsilon\xi\iota\varsigma$ or
\textit{methexis}) represent the participation relation so

\begin{center}
"$x\ \mu\ u_{F}$" reads as "$x$ participates in $u_{F}$".
\end{center}

Given a relation $\mu$, an entity $u_{F}$ is said to be a \textit{universal}
for the property $F$ (with respect to $\mu$) if it satisfies the following
universality condition:

\begin{center}
for any $x$, $x\ \mu\ u_{F}$ if and only if $F(x)$.
\end{center}

A universal representing a property should be in some sense unique. Hence
there should be an equivalence relation ($\approx$) so that universals satisfy
a uniqueness condition:

\begin{center}
if $u_{F}$ and $u_{F}^{\prime}$ are universals for the same $F$, then $u_{F}$
$\approx$ $u_{F}^{\prime}$.
\end{center}

The two criteria for a \textit{theory of universals} is that it contains a
binary relation $\mu$ and an equivalence relation $\approx$ so that with
certain properties $F$ there are associated entities $u_{F}$ satisfying the
following conditions:

(1) \textit{Universality condition}: for any $x$, $x\ \mu\ u_{F}$ iff $F(x)$, and

(2) \textit{Uniqueness condition}: if $u_{F}$ and $u_{F}^{\prime}$ are
universals for the same $F$ [i.e., satisfy (1)], then $u_{F}$ $\approx$
$u_{F}^{\prime}$.

A universal $u_{F}$ is said to be \textit{non-self-predicative} if it does not
participate in itself, i.e., $\lnot(u_{F}\
\mu
\ u_{F})$. A universal $u_{F}$ is \textit{self-predicative }if it participates
in itself, i.e., $u_{F}\
\mu
\ u_{F}$.\footnote{A self-predicative universal for some property is thus an
impredicative definition of having that property. See \cite[p. 245]%
{louie-poli:hspread} where a supremum or least upper bound is referred to as
giving an impredicative definition of being an upper bound of a subset of a
partial order. Also Michael Makkai \cite{makkai:struct} makes a similiar
remark about the universal mapping property of the natural number system.} For
the sets in an iterative set theory (Boolos \cite{boolos:iterative}), set
membership is the participation relation, set equality is the equivalence
relation, and those sets are never-self-predicative (since the set of
instances of a property is always of higher type or rank than the instances).
The universals of category theory form the "other bookend" as
always-self-predicative universals. The set-theoretical paradoxes arose from
trying to have \textit{one} theory of universals ("Frege's Paradise") where
the universals could be \textit{either }self-predicative or
non-self-predicative,\footnote{Then the universal for all the
non-self-predicative universals would give rise to Russell's Paradox since it
could not be self-predicative or non-self-predicative (Russell \cite[p.
80]{russell:prinmath}).} instead of having two opposite "bookend" theories,
one for never-self-predicative universals (set theory) and one for always
always-self-predicative universals (category theory).

For the self-predicative universals of category theory (see
\cite{maclane:algebra} or \cite{maclane:cwm} for introductions), the
participation relation is the \textit{uniquely-factors-through} relation. It
can always be formulated in a suitable category as:

\begin{center}
"$x\ \mu\ u_{F}$" means "there exists a unique arrow $x\Rightarrow u_{F}$".
\end{center}

Then $x$ is said to \textit{uniquely factor through} $u_{F}$, and the arrow
$x\Rightarrow u_{F}$ is the unique factor or participation morphism. In the
universality condition,

\begin{center}
for any $x$, $x\ \mu\ u_{F}$ if and only if $F(x)$,
\end{center}

\noindent the existence of the identity arrow $1_{u_{F}}:u_{F}\Rightarrow
u_{F}$ is the self-participation of the self-predicative universal that
corresponds with $F(u_{F})$, the self-predication of the property to $u_{F}$.
In category theory, the equivalence relation used in the uniqueness condition
is the isomorphism ($\cong$).

\section{The Hom-set definition of an adjunction}

We will later use a specific heterodox treatment of adjunctions, first
developed by Pareigis \cite{pareigis:cats-functors} and later rediscovered and
developed by Ellerman (\cite{ellerman:whatisct}, \cite{ellerman:axiomathes}),
which shows that adjoints arise by gluing together in a certain way two
universals (left and right representations). But for illustration, we start
with the standard Hom-set definition of an adjunction.

The category $Sets$ has all sets as objects and all functions between sets as
the homomorphisms so for sets $a$ and $a^{\prime}$, $\operatorname*{Hom}%
\left(  a,a^{\prime}\right)  $ is the set of functions $a\rightarrow
a^{\prime}$. In the product category $Sets\times Sets$, the objects are
ordered pairs of sets $\left(  a,b\right)  $ and homomorphism $\left(
a,b\right)  \rightarrow\left(  a^{\prime},b^{\prime}\right)  $ is just a pair
of functions $\left(  f,g\right)  $ where $f:a\rightarrow a^{\prime}$ and
$g:b\rightarrow b^{\prime}$.

For an example of an adjunction, consider the \textit{product functor}
$\times:Sets\times Sets\rightarrow Sets$ which takes a pair of sets $\left(
a,b\right)  $ to their Cartesian product $a\times b$ (set of ordered pairs of
elements from $a$ and $b$) and takes a homomorphism $\left(  f,g\right)
:\left(  a,b\right)  \rightarrow\left(  a^{\prime},b^{\prime}\right)  $ to
$f\times g:a\times b\rightarrow a^{\prime}\times b^{\prime}$ where for $x\in
a$ and $y\in b$, $f\times g:\left(  x,y\right)  \longmapsto\left(  f\left(
x\right)  ,g\left(  y\right)  \right)  $.

The maps $f:a\rightarrow a^{\prime}$ in $Sets$ go from one set to one set and
the maps $\left(  f,g\right)  :\left(  a,b\right)  \rightarrow\left(
a^{\prime},b^{\prime}\right)  $ in $Sets\times Sets$ go from a pair of sets to
a pair of sets. There is also the idea of a \textit{cone} $\left[  f,g\right]
:c\rightarrow\left(  a,b\right)  $ of maps that is a pair of maps
$f:c\rightarrow a$ and $g:c\rightarrow b$ going from one set $c$ (the point of
the cone) in $Sets$ to a pair of sets $\left(  a,b\right)  $ (the base of the
cone) in $Sets\times Sets$. Before the notion of a adjunction was defined by
Kan \cite{kan:adjoints}, the product of sets $a\times b$ was defined by its
universal mapping property. The projection maps $\pi_{a}:a\times b\rightarrow
a$ and $\pi_{b}:a\times b\rightarrow b$ define a canonical cone $\left[
\pi_{a},\pi_{b}\right]  :a\times b\rightarrow(a,b)$ that is universal in the
following sense. Given any other cone $\left[  f,g\right]  :c\rightarrow(a,b)$
from any set $c$ to $\left(  a,b\right)  $, there is a unique homomorphism
$\left\langle f,g\right\rangle :c\rightarrow a\times b$ in $Sets$ such that
the two triangles in the following diagram commute.%

\begin{center}
\includegraphics[
height=1.695in,
width=3.1142in
]%
{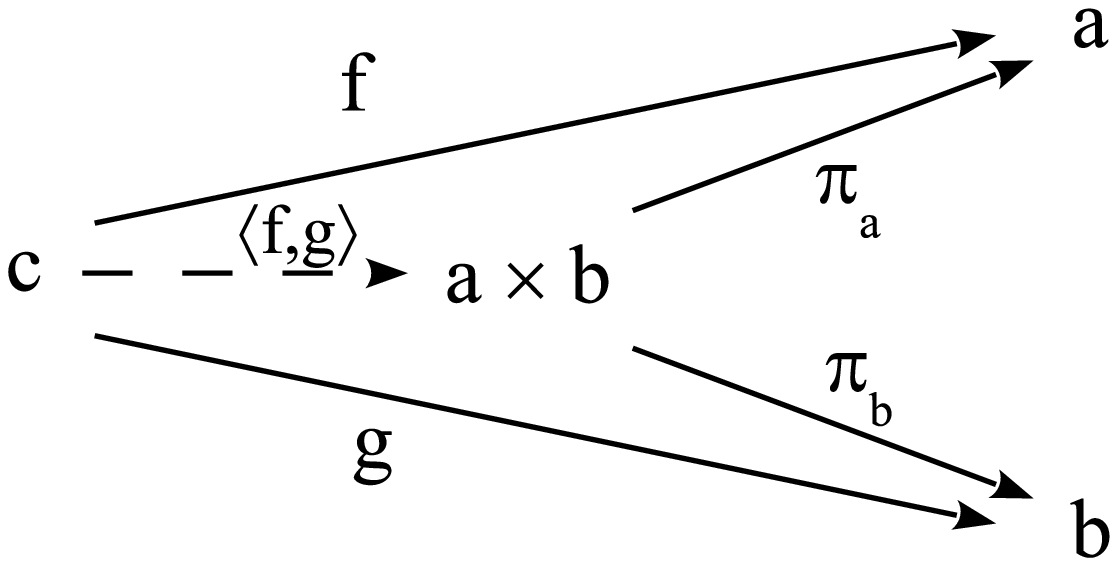}%
\end{center}

\begin{center}
Figure 1: Universal mapping property for the product of sets.
\end{center}

In terms of the self-predicative universals considered in the last section,
the property in question is the property of being a cone $\left[  f,g\right]
:c\rightarrow(a,b)$ to $\left(  a,b\right)  $ from any set $c$. The canonical
cone of projections $\left[  \pi_{a},\pi_{b}\right]  :a\times b\rightarrow
(a,b)$ is the self-predicative universal for that property. The participation
relation $\left[  f,g\right]  \ \mu\ \left[  \pi_{a},\pi_{b}\right]  $ is
defined as "uniquely factoring through" (as in Figure 1). The universal
mapping property of the product can then be restated as the universality
condition: For any cone $\left[  f,g\right]  $ from any set to a pair of sets,

\begin{center}
$\left[  f,g\right]  \ \mu\ \left[  \pi_{a},\pi_{b}\right]  $ if and only if
$\left[  f,g\right]  $ is a cone to $\left(  a,b\right)  $.

UMP of $a\times b$ stated as a universality condition.
\end{center}

The Hom-set definition of the adjunction for the product functor uses the
auxiliary device of a diagonal functor to avoid mentioning the cones and to
restrict attention only to the Hom-sets of the two categories. The
\textit{diagonal functor} $\Delta:Sets\rightarrow Sets\times Sets$ in the
opposite direction of the product functor just doubles everything so
$\Delta\left(  c\right)  =\left(  c,c\right)  $ and $\Delta\left(  f\right)
=\left(  f,f\right)  $. Then the product functor is said to be the
\textit{right adjoint} of the diagonal functor, the diagonal functor is said
to be the\textit{ left adjoint} of the product functor, and the two functors
together form an \textit{adjunction} if there is a natural isomorphism between
the Hom-sets as follows:

\begin{center}
$\operatorname*{Hom}_{Sets\times Sets}\left(  \Delta\left(  c\right)  ,\left(
a,b\right)  \right)  \cong\operatorname*{Hom}_{Sets}\left(  c,a\times
b\right)  $.

Hom-set definition of the adjunction between the product and diagonal functors.
\end{center}

The diagonal functor $\Delta:Sets\rightarrow Sets\times Sets$ also has a
(rather trivial) UMP that can be stated in terms of cones $c\rightarrow(a,b)$
except now we fix $c$ and let $\left(  a,b\right)  $ vary. There is the
canoncial cone $\left[  1_{c},1_{c}\right]  :c\rightarrow\left(  c,c\right)  $
and it is universal in the following sense. For any cone $\left[  f,g\right]
:c\rightarrow\left(  a,b\right)  $ from the given $c$ to any pair of sets
$\left(  a,b\right)  $, there is a unique homomorphism in $Sets\times Sets$,
namely $\left(  f,g\right)  :\left(  c,c\right)  \rightarrow\left(
a,b\right)  $ that factors through the canonical cone $c\rightarrow\left(
c,c\right)  $.%

\begin{center}
\includegraphics[
height=1.3015in,
width=3.4177in
]%
{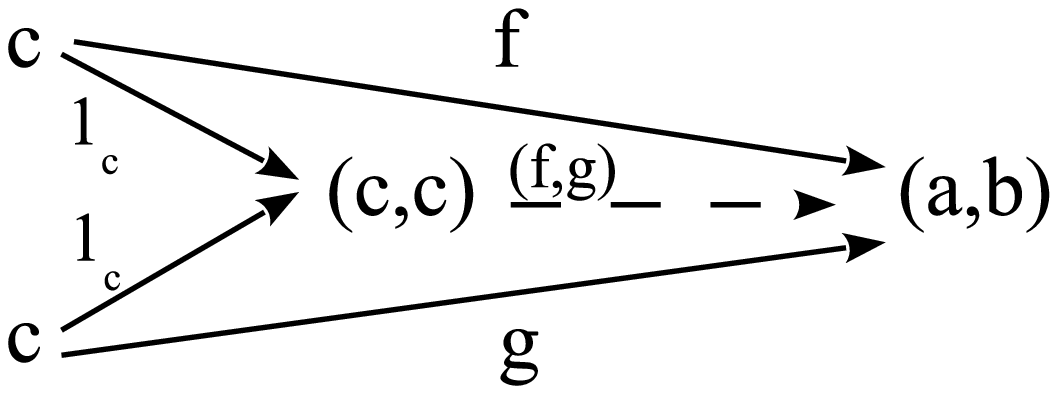}%
\end{center}

\begin{center}
Figure 2: Universal mapping property for diagonal functor.
\end{center}

This product-diagonal adjunction illustrates the general $\operatorname*{Hom}%
$-set definition. Given functors $F:\mathbb{X}\rightarrow\mathbb{A}$ and
$G:\mathbb{A}\rightarrow\mathbb{X}$ going each way between categories
$\mathbb{X}$ and $\mathbb{A}$, they form an adjunction if there is a natural
isomorphism (for objects $X\in\mathbb{X}$ and $A\in\mathbb{A}$):

\begin{center}
$\operatorname*{Hom}_{\mathbb{A}}\left(  F\left(  X\right)  ,A\right)
\cong\operatorname*{Hom}_{\mathbb{X}}\left(  X,G\left(  A\right)  \right)  $

$\operatorname*{Hom}$-set definition of an adjunction.
\end{center}

To further analyze adjoints, we need the notion of a "heteromorphism."

\section{Heteromorphisms and adjunctions}

We have seen that there are two UMPs (often one is trivial like $\Delta(c)$ in
the above example) involved in an adjunction and that the object-to-object
maps were always within one category, e.g., in the "Hom-sets" of one category
or the other. Using object-to-object maps between objects of
\textit{different} categories (properly called "heteromorphisms" or "chimera
morphisms"), the notion of an adjunction can be factored into two
representations (or "half-adjunctions" in Ellerman \cite[p. 158]%
{ellerman:whatisct}), each of which expresses a universal mapping property.

We have already seen one standard example of a \textit{heteromorphism} or
\textit{het}, namely a cone $\left[  f,g\right]  :c\rightarrow\left(
a,b\right)  $ that goes from an object in $Sets$ to an object in $Sets\times
Sets$. The hets are contrasted with the homs or homomorphisms between objects
in the same category. To keep them separate in our notation, we will
henceforth use single arrows $\longrightarrow$ for hets and double arrows
$\Rightarrow$ for homs.\footnote{The hets between objects of different
categories are represented as single arrows ($\rightarrow$) while the
homomorphisms or homs between objects in the same category are represented by
double arrows ($\Rightarrow$). The functors between whole categories are also
represented by single arrows ($\rightarrow$). One must be careful not to
confuse a functor $F:\mathbb{X}\rightarrow\mathbb{A}$ from a category
$\mathbb{X}$ to a category $\mathbb{A}$ with its action on an object
$X\in\mathbb{X}$ which would be symbolized $X\longmapsto F(X)$. Moreover since
a functor often has a canonical definition, there may well be a canonical het
$X\rightarrow F(X)$ or $X\leftarrow F\left(  X\right)  $ but such hets are no
part of the definition of the functor itself.} Then the UMP for the product
functor can be represented as follows.

\begin{center}%
\begin{center}
\includegraphics[
height=1.5091in,
width=2.6515in
]%
{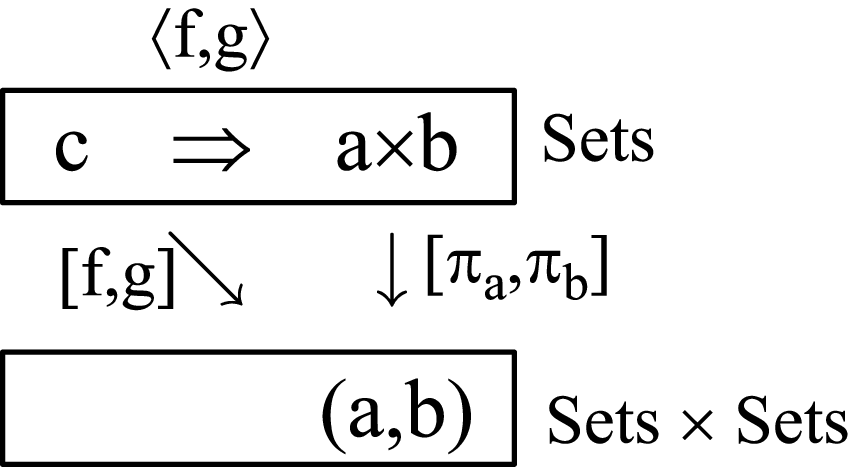}%
\end{center}

Figure 3: UMP for the product functor
\end{center}

\noindent It should be particularly noted that this het-formulation of the UMP
for the product does not involve the diagonal functor. If we associate with
each $c\in Sets$ and each $\left(  a,b\right)  \in Sets\times Sets$, the set
$\operatorname*{Het}\left(  c,\left(  a,b\right)  \right)  $ of cones or hets
$\left[  f,g\right]  :c\rightarrow\left(  a,b\right)  $ then this defines a
$\operatorname*{Het}$\textit{-bifunctor} in the same manner as the usual
$\operatorname*{Hom}$-bifunctor $\operatorname*{Hom}_{Sets}\left(
a,a^{\prime}\right)  $ or $\operatorname*{Hom}_{Sets\times Sets}\left(
\left(  a,b\right)  ,\left(  a^{\prime},b^{\prime}\right)  \right)  $ [see the
appendix for more details]. Then the UMP for the product functor gives a
natural isomorphism based on the pairing: $\left[  f,g\right]  \mapsto
\left\langle f,g\right\rangle $, so that the $Sets$-valued functor
$\operatorname*{Het}\left(  c,\left(  a,b\right)  \right)  $ is said to be
\textit{represented on the right} by the $Sets$-valued $\operatorname*{Hom}%
_{Sets}\left(  c,a\times b\right)  $:

\begin{center}
$\operatorname*{Het}\left(  c,\left(  a,b\right)  \right)  \cong%
\operatorname*{Hom}_{Sets}\left(  c,a\times b\right)  $

Right representation of the hets $c\rightarrow\left(  a,b\right)  $ by the
homs $c\Rightarrow a\times b$.
\end{center}

The trivial UMP for the diagonal functor can also be stated in terms of the
cone-hets without reference to the product functor.

\begin{center}%
\begin{center}
\includegraphics[
height=1.3067in,
width=2.7951in
]%
{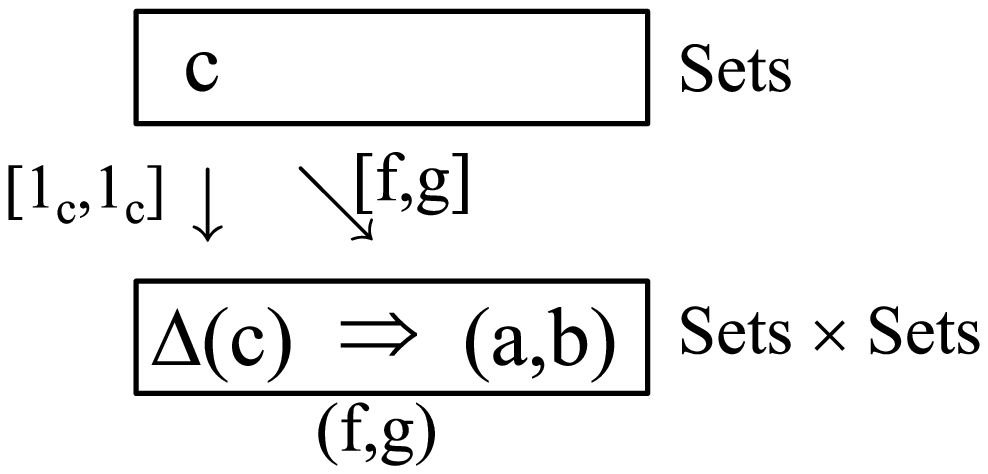}%
\end{center}

Figure 4: UMP for the diagonal functor
\end{center}

\noindent This UMP for the diagonal functor gives a natural isomorphism based
on the pairing $\left(  f,g\right)  \longmapsto\left[  f,g\right]  $, so the
$Sets$-valued functor $\operatorname*{Het}\left(  c,\left(  a,b\right)
\right)  $ is said to be \textit{represented on the left} by the $Sets$-valued
$\operatorname*{Hom}_{Sets\times Sets}\left(  \left(  c,c\right)  ,\left(
a,b\right)  \right)  $:

\begin{center}
$\operatorname*{Hom}_{Sets\times Sets}\left(  \left(  c,c\right)  ,\left(
a,b\right)  \right)  \cong$ $\operatorname*{Het}\left(  c,\left(  a,b\right)
\right)  $

Left representation of the hets $c\rightarrow\left(  a,b\right)  $ by the homs
$\left(  c,c\right)  \Rightarrow\left(  a,b\right)  $.
\end{center}

Then the right and left representations of the hets $\operatorname*{Het}%
\left(  c,\left(  a,b\right)  \right)  $ can be combined to obtain as a
consequence the Hom-set definition of the adjunction between the product and
diagonal functors:

\begin{center}
$\operatorname*{Hom}_{Sets\times Sets}\left(  \left(  c,c\right)  ,\left(
a,b\right)  \right)  \cong$ $\operatorname*{Het}\left(  c,\left(  a,b\right)
\right)  \cong\operatorname*{Hom}_{Sets}\left(  c,a\times b\right)  $

Heteromorphic presentation of the product-diagonal adjunction.
\end{center}

In the general case of adjoint functors $F:\mathbb{X}\rightleftarrows
\mathbb{A}:G$, the hets $\operatorname*{Het}\left(  X,A\right)  $ from objects
$X\in\mathbb{X}$ to objects $A\in\mathbb{A}$ have left and right representations:

\begin{center}
$\operatorname*{Hom}_{\mathbb{A}}\left(  F\left(  X\right)  ,A\right)
\cong\operatorname*{Het}\left(  X,A\right)  \cong\operatorname*{Hom}%
_{\mathbb{X}}\left(  X,G\left(  A\right)  \right)  $

Heteromorphic presentation of a general adjunction.
\end{center}

This is the heterodox treatment of an adjunction first published by Pareigis
\cite[pp. 60-1]{pareigis:cats-functors} and later rediscovered and developed
by the author (\cite[p. 130]{ellerman:whatisct} and \cite{ellerman:axiomathes}%
). It is "heterodox" since the morphisms between the objects of
\textit{different} categories are not "officially" recognized in the standard
presentations of category theory (e.g., \cite{maclane:cwm} or
\cite{awodey:cat}) even though such hets are a common part of mathematical
practice (see the appendix for further discussion). Hence the standard Hom-set
definition of an adjunction just deletes the $\operatorname*{Het}$-middle-term
$\operatorname*{Het}\left(  X,A\right)  $ to obtain just the het-free or
homs-only presentation of an adjunction.

The important advance of the heteromorphic treatment of an adjunction is that
the adjunction can be parsed or factored into two parts, the left and right
representations, each of which only involves one of the $\operatorname*{Hom}$-functors.

\begin{center}%
\begin{center}
\includegraphics[
height=1.5203in,
width=4.126in
]%
{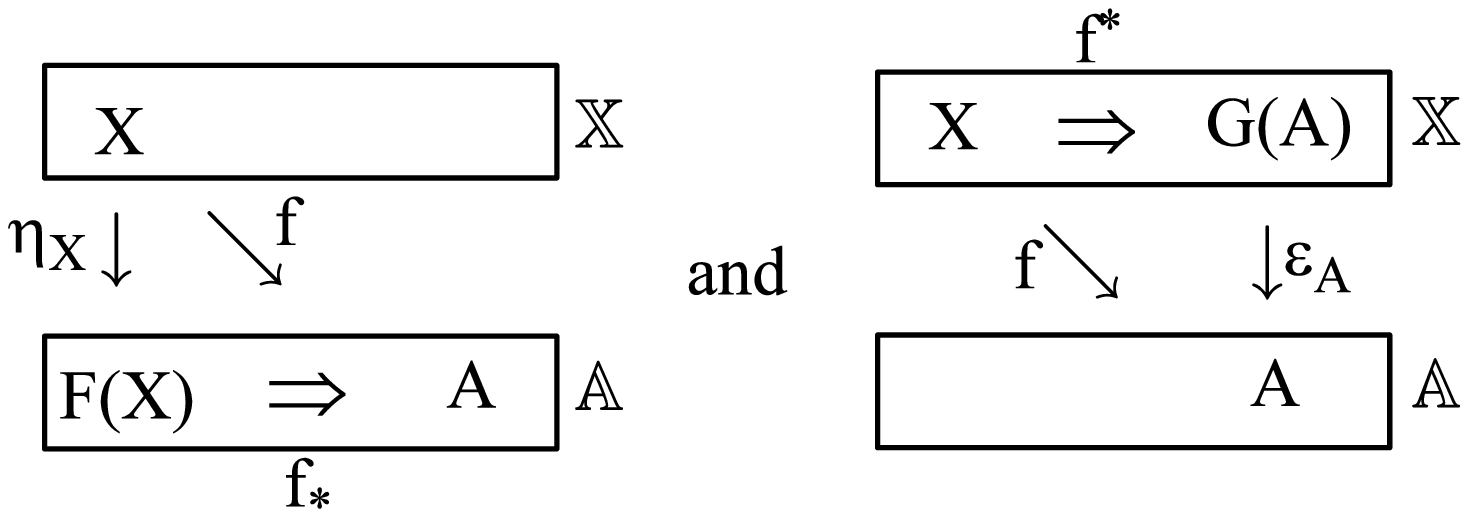}%
\end{center}

Figure 5: Left and right representations each involving on one of the adjoints
$F$ and $G$.
\end{center}

\noindent Moreover, the diagrams for the two representations can be glued
together at the diagonal het $\searrow^{f}$ into one diagram to give the
simple \textit{adjunctive square diagram} for an adjunction.

\begin{center}%
\begin{center}
\includegraphics[
height=1.7841in,
width=2.0548in
]%
{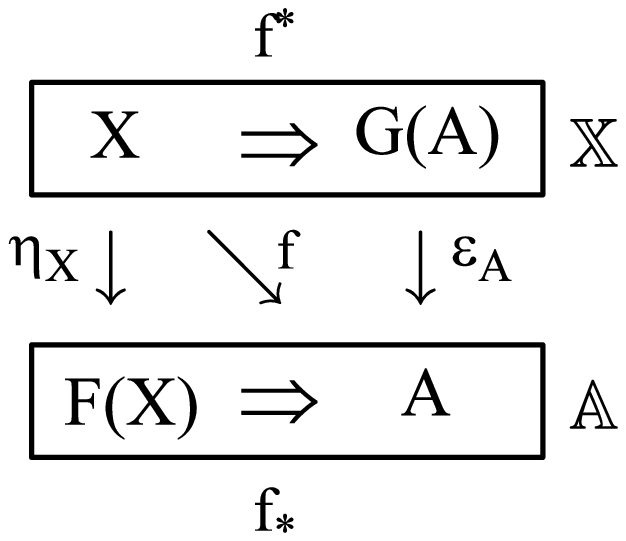}%
\end{center}

Figure 6: Adjunctive square diagram for the het-treatment of an adjunction.
\end{center}

Every adjunction can be represented (up to isomorphism) in this manner
\cite[p. 147]{ellerman:whatisct} so the molecule of an adjunction can be split
into two atoms, each of which is a (left or right) representation of a
$\operatorname*{Het}$-functor. This means that the importance and ubiquity of
adjunctions (emphasized above) also passes to the atoms, left or right
representations, that make up those molecules. Moreover, it should be noted
that each left or right representation defines a self-predicative universal as
indicated in the previous example of the het-cones $c\rightarrow\left(
a,b\right)  $.

The main point of this paper is that those atoms, the left and right
representations can be recombined in a new way to define a "recombinant
construction" cognate to an adjunction, and that is the concept of a "brain functor."

\section{Brain functors}

In many adjunctions, the important fact is expressed by either the left or
right representation (e.g., the UMP for the product functor or for the
free-group functor considered in the Appendix), with no need for the
"auxiliary device" (such as a diagonal or forgetful functor) of the other
representation used to express the adjunction in a het-free manner.

Another payoff from analyzing the important but molecular concept of an
adjunction into two atomic representations is that we can then reassemble
those atomic parts in a new way to define the cognate concept speculatively
named a "brain functor."

The basic intuition is to think of one category $\mathbb{X}$ in a
representation as the "environment" and the other category $\mathbb{A}$ as an
"organism." Instead of representations within \textit{each} category of the
hets going \textit{one way} between the categories (as in an adjunction),
suppose the hets going \textit{both ways} were represented within \textit{one}
of the categories (the "organism").

Intuitively, a het from the environment to the organism is say, a visual or
auditory stimulus. Then a left representation would play the role of the brain
in providing the re-cognition or perception (expressed by the
intentionality-of-perception slogan: "seeing is seeing-as") of the stimulus as
a perception of, say, a tree where the internal re-cognition is represented by
the homomorphism $\Rightarrow$ inside the "organism" category.%

\begin{center}
\includegraphics[
height=1.7789in,
width=3.4636in
]%
{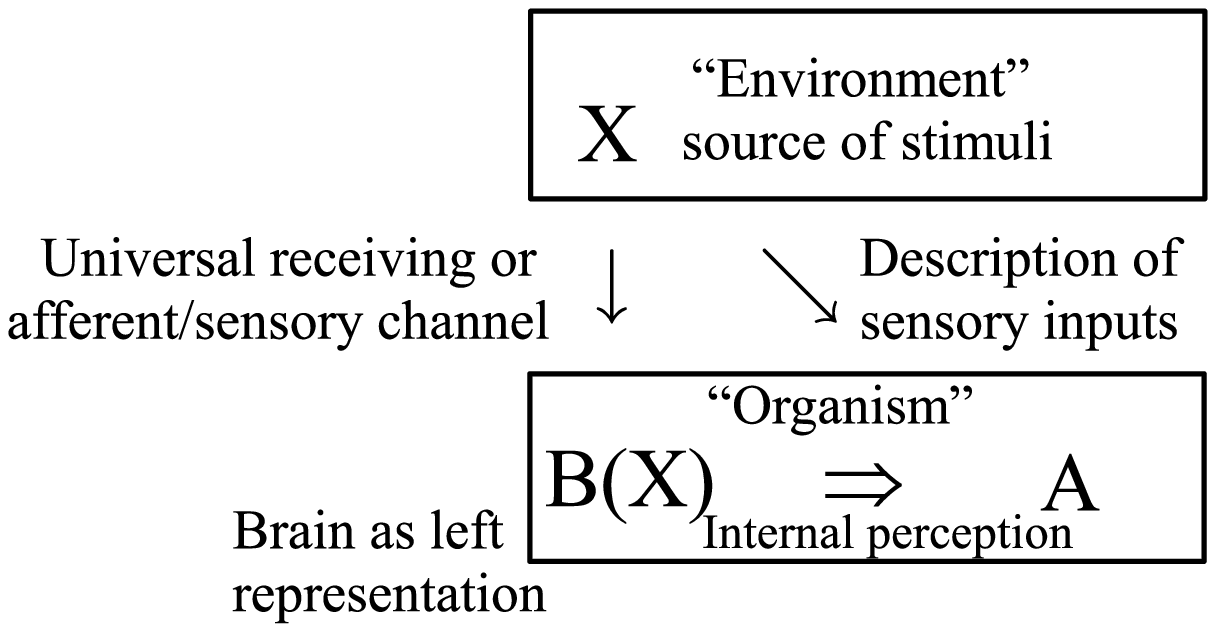}%
\end{center}

\begin{center}
Figure 7: Perceiving brain presented as a left representation.
\end{center}

Perhaps not surprisingly, this mathematically models the old philosophical
theme in the Platonic tradition that external stimuli do not give knowledge;
the stimuli only trigger the internal perception, recognition, or recollection
(as in Plato's \textit{Meno}) that is knowledge. In \textit{De Magistro} (The
Teacher), the neo-Platonic Christian philosopher Augustine of Hippo developed
an argument (in the form of a dialogue with his son Adeodatus) that as
teachers teach, it is only the student's internal appropriation of what is
taught that gives understanding.

\begin{quotation}
\noindent Then those who are called pupils consider within themselves whether
what has been explained has been said truly; looking of course to that
interior truth, according to the measure of which each is able. Thus they
learn,\ldots. But men are mistaken, so that they call those teachers who are
not, merely because for the most part there is no delay between the time of
speaking and the time of cognition. And since after the speaker has reminded
them, the pupils quickly learn within, they think that they have been taught
outwardly by him who prompts them. (Augustine, \textit{De Magistro}, Chapter XIV)
\end{quotation}

The basic point is the active role of the mind in generating understanding
(represented by the internal hom). This is clear even at the simple level of
understanding spoken words. We hear the auditory sense data of words in a
completely strange language as well as the words in our native language. But
the strange words bounce off our minds, like @\#\$\%\symbol{94}, with no
resultant understanding while the words in a familiar language prompt an
internal process of generating a meaning so that we understand the words. Thus
it could be said that "understanding a language" means there is a left
representation for the heard statements in that language, but there is no such
internal re-cognition mechanism for the heard auditory inputs in a strange language.

Dually, there are also hets going the other way from the "organism" to the
"environment" and there is a similar distinction between mere behavior (e.g.,
a reflex) and an action that expresses an intention. Mathematically that is
described by dualizing or turning the arrows around which gives an acting
brain presented as a right representation.%

\begin{center}
\includegraphics[
height=1.6994in,
width=3.2932in
]%
{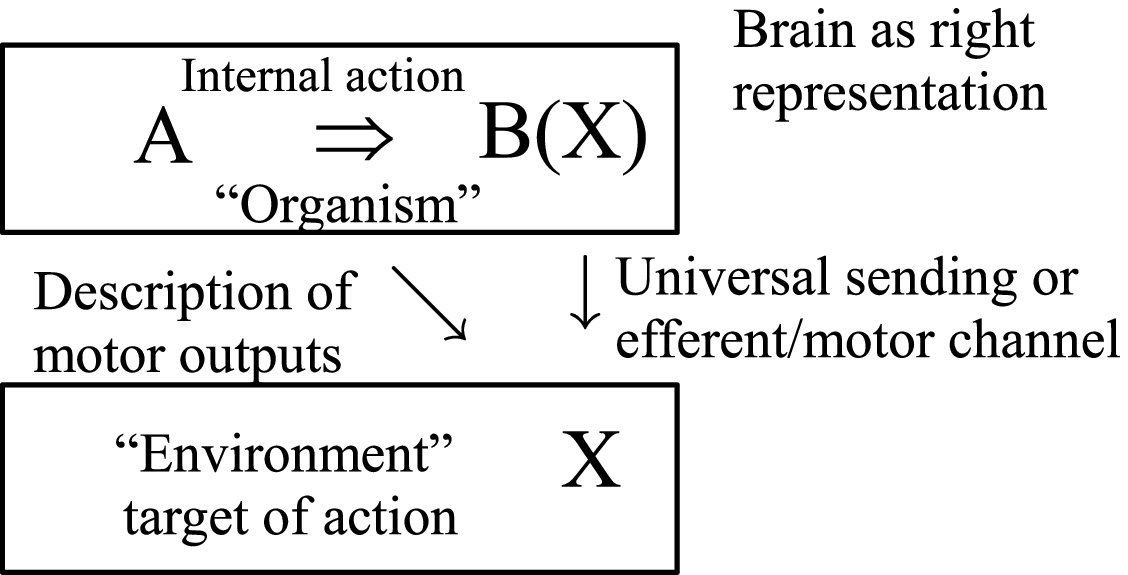}%
\end{center}

\begin{center}
Figure 8: Acting brain as a right representation.
\end{center}

In the heteromorphic treatment of adjunctions, an adjunction arises when the
hets from one category $\mathbb{X}$ to another category $\mathbb{A}$,
$\operatorname*{Het}(X,A)$ for $X\in\mathbb{X}$ and $A\in\mathbb{A}$, have a
right representation, $\operatorname*{Het}(X,A)\cong\operatorname*{Hom}%
_{\mathbb{X}}(X,G(A))$, and a left representation, $\operatorname*{Hom}%
_{\mathbb{A}}(F(X),A)\cong\operatorname*{Het}(X,A)$. But instead of taking the
same set of hets as being represented by two different functors on the right
and left, suppose we consider a single functor $B(X)$ that represents the hets
$\operatorname*{Het}(X,A)$ on the left:

\begin{center}
$\operatorname*{Het}(X,A)\cong\operatorname*{Hom}_{\mathbb{A}}(B(X),A)$,
\end{center}

\noindent and represents the hets $\operatorname*{Het}(A,X)$ [going in the
opposite direction] on the right:

\begin{center}
$\operatorname*{Hom}_{\mathbb{A}}(A,B(X))\cong\operatorname*{Het}(A,X)$.
\end{center}

If the hets each way between two categories are represented by the same
functor $B(X)$ as left and right representations, then that functor is said to
be a \textit{brain functor}. Thus instead of a pair of functors being adjoint,
we have a single functor $B(X)$ with values within one of the categories (the
"organism") as representing the two-way interactions, "perception" and
"action," between that category and another one (the "environment"). The use
of the adjective "brain" is quite deliberate (as opposed to say "mind") since
the universal hets going each way between the "organism" and "environment" are
part of the definition of left and right representations. In particular, it
should be noted how the "turn-around-the-arrows" category-theoretic duality
provides a mathematical model for the type of "duality" between:

\begin{itemize}
\item sensory or afferent systems (brain furnishing the left representation of
the environment to organism heteromorphisms), and

\item motor or efferent systems (brain furnishing the right representation of
the organism to environment heteromorphisms).
\end{itemize}

In view of this application, those two universal hets, representing the
afferent and efferent nervous systems, might be denoted $\alpha_{X}$ and
$\varepsilon_{X}$ as in the following diagrams for the two representations.

\begin{center}%
\begin{center}
\includegraphics[
height=1.6129in,
width=4.3267in
]%
{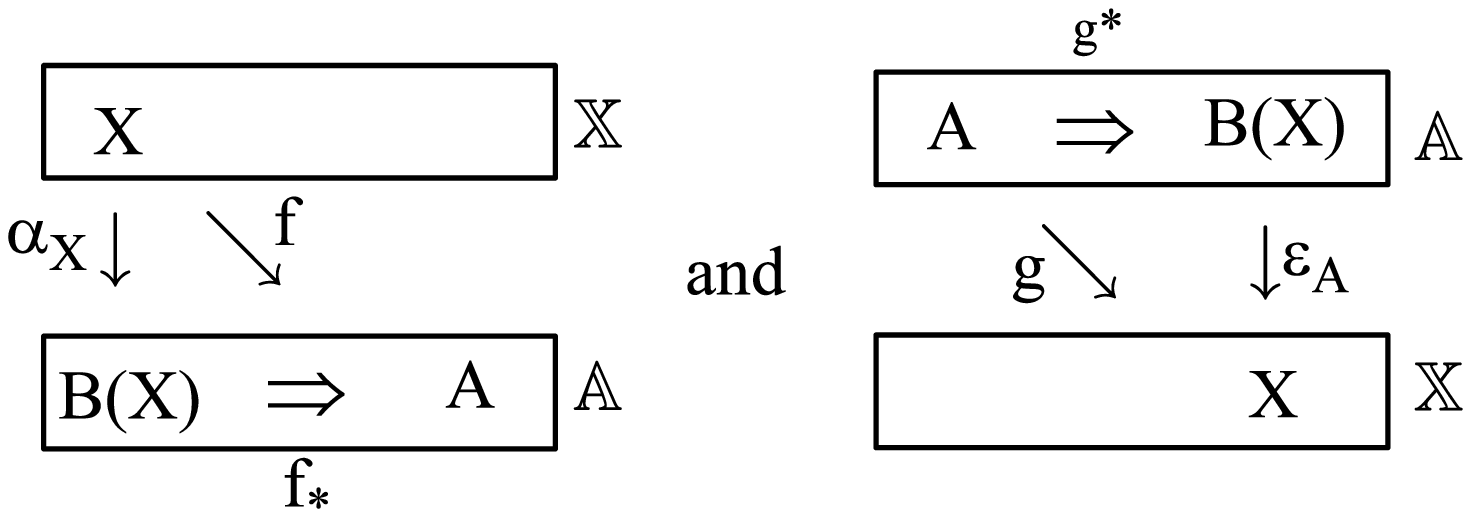}%
\end{center}

Figure 9: Left and right representation diagrams for the brain functor
$B:\mathbb{X}\rightarrow\mathbb{A}$.
\end{center}

We have seen how the adjunctive square diagram for an adjunction can be
obtained by gluing together the left and right representation diagrams at the
common diagonal $\searrow^{f}$. The diagram for a brain functor is obtained by
gluing together the diagrams for the left and right representations at the
common values of the brain functor $B(X)$. If we think of the diagram for a
representation as right triangle, then the adjunctive square diagram is
obtained by gluing two triangles together on the hypotenuses, and the diagram
for the brain functor is obtained by gluing two triangles together at the
right angle vertices to form the \textit{butterfly diagram. }

\begin{center}%
\begin{center}
\includegraphics[
height=2.1119in,
width=2.9205in
]%
{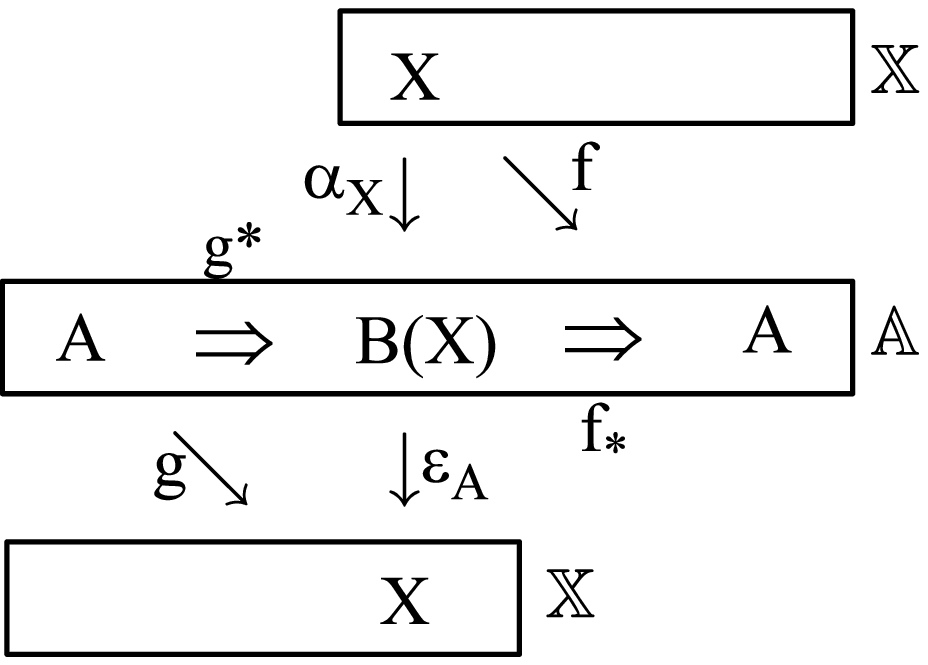}%
\end{center}

Figure 10: Butterfly diagram combining two representations at the common
$B\left(  X\right)  $
\end{center}

If both the triangular "wings" could be filled-out as adjunctive squares, then
the brain functor would have left and right adjoints. Thus all functors with
both left and right adjoints are brain functors (although not vice-versa). The
previous example of the diagonal functor $\Delta:Sets\rightarrow Sets\times
Sets$ is a brain functor since the product functor $\times\left(  a,b\right)
=a\times b$ is the right adjoint, and the coproduct or disjoint union functor
$%
{\textstyle\biguplus}
\left(  a,b\right)  =a%
{\textstyle\biguplus}
b$ is the left adjoint. The underlying set functor (see Appendix) that takes a
group $G$ to its underlying set $U\left(  G\right)  $ is a rather trivial
example of a brain functor that does not arise from having both a left and
right adjoint. It has a left adjoint (the free group functor) so $U$ provides
a right representation for the set-to-group maps or hets $X\rightarrow G$.
Also it trivially provides a left representation for the hets $G\rightarrow X$
but has no right adjoint.

In the butterfly diagram below, we have labelled the diagram for the brain as
the language faculty for understanding and producing speech.
\begin{center}
\includegraphics[
height=2.1741in,
width=3.6348in
]%
{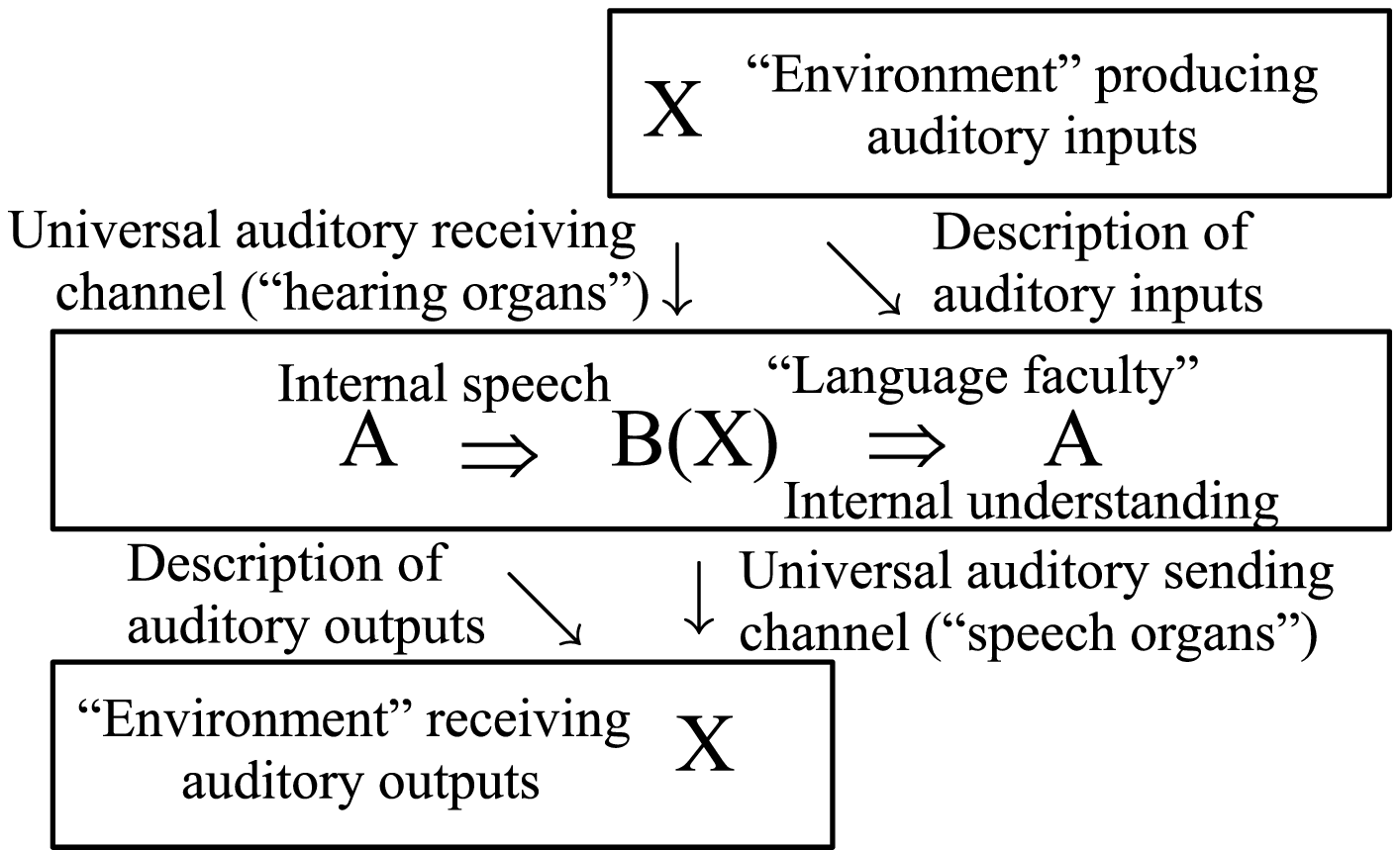}%
\end{center}

\begin{center}
Figure 11: Brain functor butterfly diagram interpreted as language faculty.
\end{center}

Wilhelm von Humboldt recognized the symmetry between the speaker and listener,
which in the same person is abstractly represented as the dual functions of
the "selfsame power" of the language faculty in the above butterfly diagram.

\begin{quotation}
\noindent Nothing can be present in the mind (Seele) that has not originated
from one's own activity. Moreover understanding and speaking are but different
effects of the selfsame power of speech. Speaking is never comparable to the
transmission of mere matter (Stoff). In the person comprehending as well as in
the speaker, the subject matter must be developed by the individual's own
innate power. What the listener receives is merely the harmonious vocal
stimulus.\cite[p. 102]{humboldt:conformation}
\end{quotation}

\section{A mathematical example of a brain functor}

A non-trivial mathematical example of a brain functor is provided by the
functor taking a finite set of vector spaces $\{V_{i}\}_{i=1,\ldots,n}$ over
the same field (or $R$-modules over a ring $R$) to the product $%
{\textstyle\prod_{i}}
V_{i}$ of the vector spaces. Such a product is also the coproduct $%
{\textstyle\sum_{i}}
V_{i}$ \cite[p. 173]{hungerford:algebra} and that space may be written as the biproduct:

\begin{center}
$V_{1}\oplus\ldots\oplus V_{n}\cong%
{\textstyle\prod\nolimits_{i}}
V_{i}\cong\sum_{i}V_{i}$.
\end{center}

The het from a set of spaces $\{V_{i}\}$ to a single space $V$ is a
\textit{cocone} of vector space maps $\{V_{i}\Rightarrow V\}$ and the
canonical such het is the set of canonical injections $\{V_{i}\Rightarrow
V_{1}\oplus\ldots\oplus V_{n}\}$ (taking the "brain" as a coproduct) with the
"brain" at the point of the cocone. The perception left representation then
might be taken as conceptually representing the function of the brain as
integrating multiple sensory inputs into an interpreted
perception.\footnote{The cocones and cones are represented in the diagrams
using cone shapes.}%

\begin{center}
\includegraphics[
height=1.9069in,
width=3.8692in
]%
{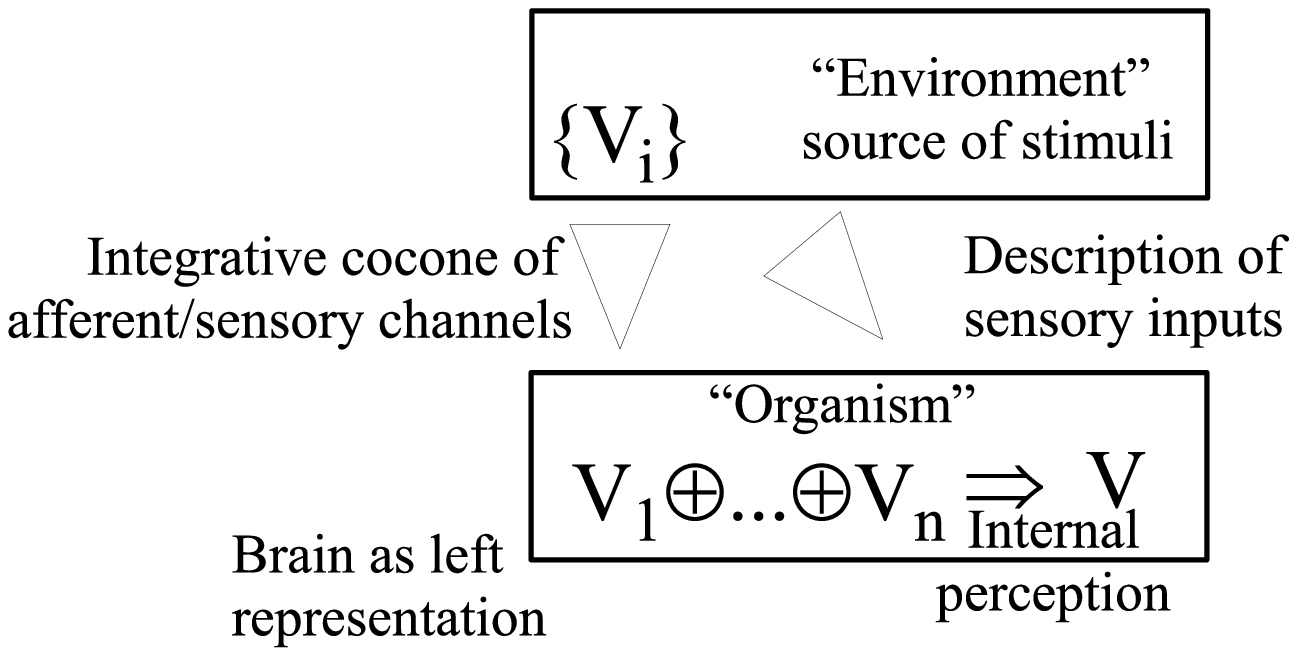}%
\end{center}

\begin{center}
Figure 12: Brain as integrating sensory inputs into a perception.
\end{center}

Dually, a het from single space $V$ to a set of vector spaces $\{V_{i}\}$ is a
\textit{cone} $\{V\Rightarrow V_{i}\}$with the single space $V$ at the point
of the cone, and the canonical het is the set of canonical projections (taking
the "brain" as a product) with the "brain" as the point of the cone:
$\{V_{1}\oplus\ldots\oplus V_{n}\Rightarrow Vi\}$. The action right
representation then might be taken as conceptually representing the function
of the brain as integrating or coordinating multiple motor outputs in the
performance of an action.%

\begin{center}
\includegraphics[
height=1.7071in,
width=3.6608in
]%
{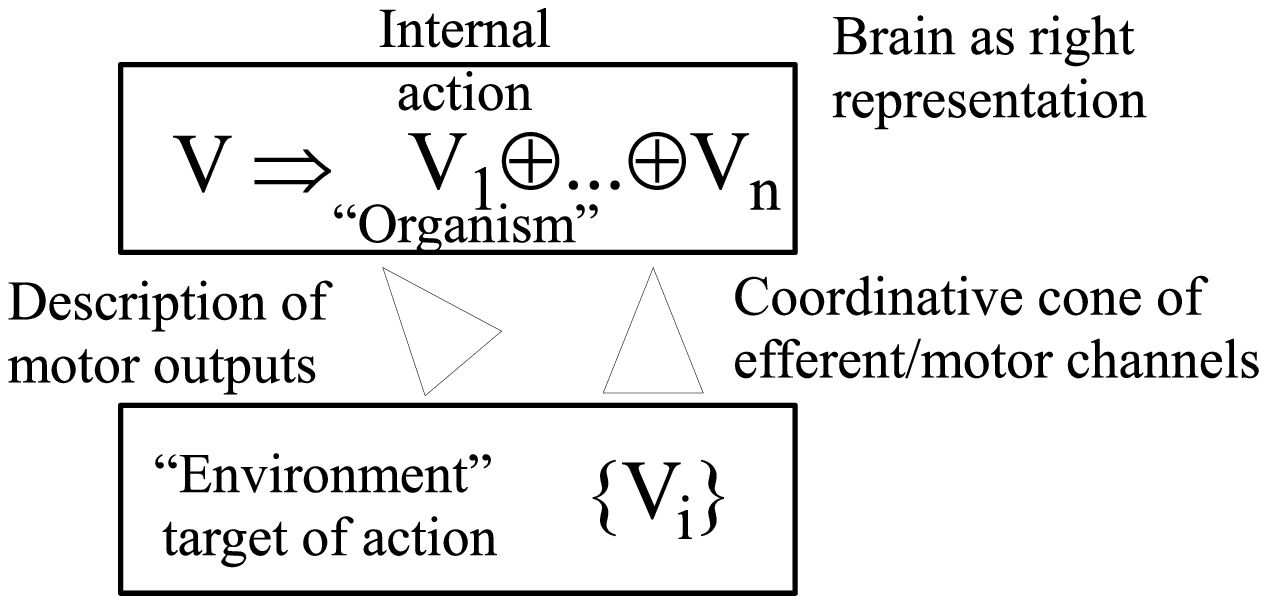}%
\end{center}

\begin{center}
Figure 13: Brain as coordinating motor outputs into an action.
\end{center}

Putting the two representations together gives the butterfly diagram for a brain.%

\begin{center}
\includegraphics[
height=2.4872in,
width=3.9115in
]%
{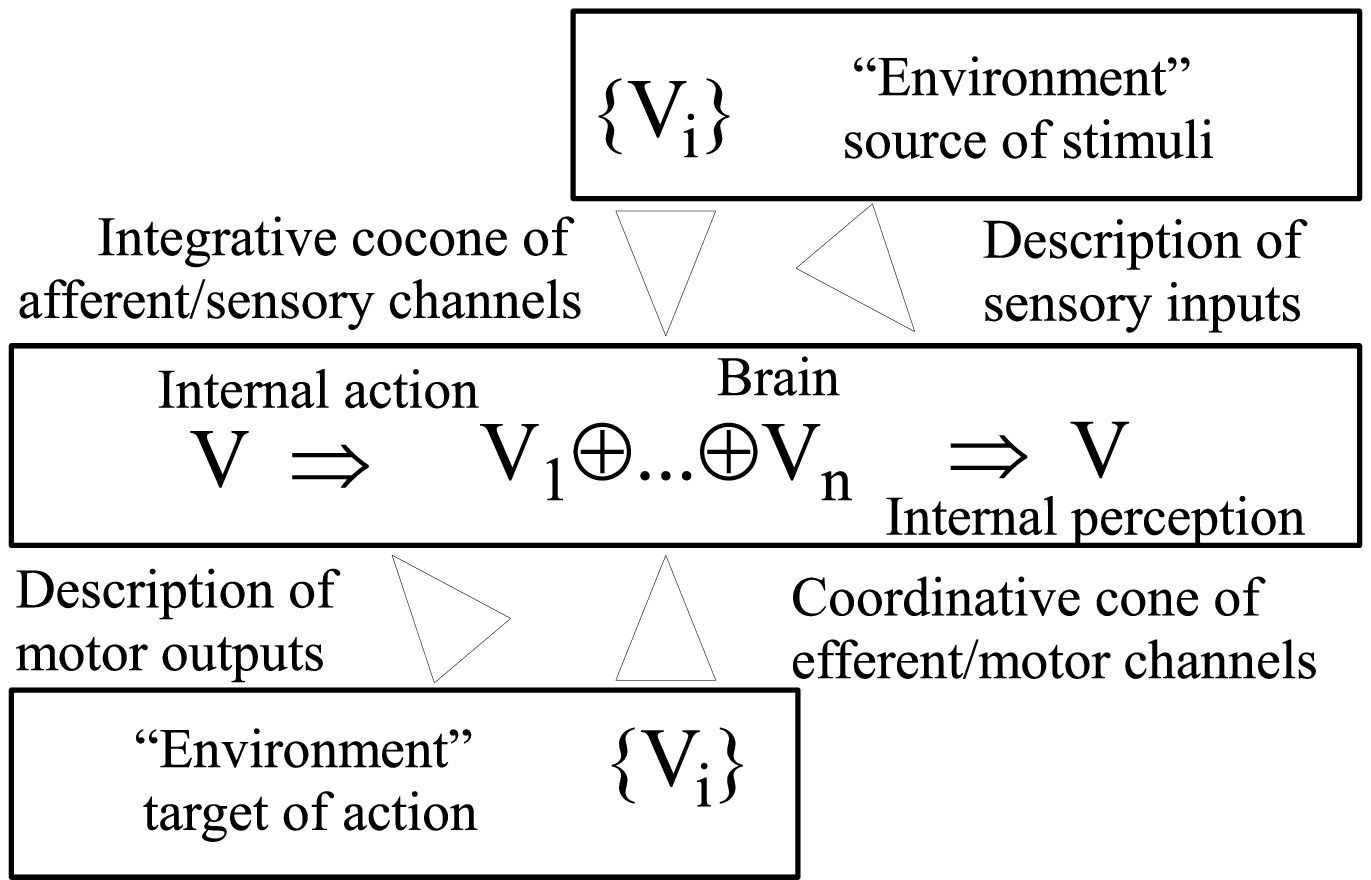}%
\end{center}

\begin{center}
Figure 14: Conceptual model of a perceiving and acting brain.
\end{center}

\noindent This gives a conceptual model of a single organ that integrates
sensory inputs into a perception and coordinates motor outputs into an action,
i.e., a brain.

\section{Conclusion}

In view of the success of category theory in modern mathematics, it is
perfectly natural to try to apply it in the life and cognitive sciences. Many
different approaches need to be tried to see which ones, if any, will find
"where theory lives" (and will be something more than just applying biological
names to bits of pure math). The approach developed here differs from other
approaches in several ways, but the most basic difference is the use of
heteromorphisms to represent interactions between quite different entities
(i.e., objects in different categories). Heteromorphisms also provide the
natural setting to formulate universal mapping problems and their solutions as
left or right representations of hets. In spite of abounding in the wilds of
mathematical practice, hets are not recognized in the orthodox presentations
of category theory. One consequence is that the notion of an adjunction
appears as one atomic concept that cannot be factored into separate parts. But
that is only a artifact of the homs-only treatment. The heteromorphic
treatment shows that an adjunction factors naturally into a left and right
representation of the hets going from one category to another--where, in
general, one representation might exist without the other. One benefit of this
heteromorphic factorization is that the two atomic concepts of left and right
representations can then be recombined in a new way to form the cognate
recombinant concept of a brain functor. The main conclusion of the paper is
that this concept of a brain functor seems to fit very well as an abstract and
conceptual but non-trivial description of the dual universal functions of a
brain, perception (using the sensory or afferent systems) and action (using
the motor or efferent systems).

\section{Mathematical Appendix: Are hets really necessary in category theory?}

Since the concept of a brain functor requires hets for its formulation, it is
important to consider the role of hets in category theory. The homomorphisms
or homs between the objects of a category $\mathbb{X}$ are given by a hom
bifunctor $\operatorname*{Hom}_{\mathbb{X}}:\mathbb{X}^{op}\times
\mathbb{X}\rightarrow Sets$. In the same manner, the heteromorphisms or hets
from the objects of a category $\mathbb{X}$ to the objects of a category
$\mathbb{A}$ are given by a het bifunctor $\operatorname*{Het}:\mathbb{X}%
^{op}\times\mathbb{A}\rightarrow Sets$.\footnote{Although often with a
somewhat different interpretation, the $Sets$-valued \textit{profunctors}
\cite{kelly:enrich}, \textit{distributors }\cite{benabou:distributors}, or
\textit{correspondences} \cite[p. 96]{lurie:highertopos} are formally the same
as het bifunctors.}

The $\operatorname*{Het}$-bifunctor gives the rigorous way to handle the
composition of a het $f:x\rightarrow a$ in $\operatorname*{Het}\left(
x,a\right)  $ [thin arrows $\rightarrow$ for hets] with a homomorphism or hom
$g:x^{\prime}\Longrightarrow x$ in $X$ [thick Arrows $\Longrightarrow$ for
homs] and a hom $h:a\Longrightarrow a^{\prime}$ in $A$. For instance, the
composition $x^{\prime}\overset{g}{\Longrightarrow}x\overset{f}{\rightarrow}a$
is the het that is the image of $f$ under the map: $\operatorname*{Het}\left(
g,a\right)  :\operatorname*{Het}\left(  x,a\right)  \rightarrow
\operatorname*{Het}\left(  x^{\prime},a\right)  $. Similarly, the composition
$x\overset{f}{\rightarrow}a\overset{h}{\Longrightarrow}a^{\prime}$ is the het
that is the image of $f$ under the map: $\operatorname*{Het}\left(
x,h\right)  :\operatorname*{Het}\left(  x,a\right)  \rightarrow
\operatorname*{Het}\left(  x,a^{\prime}\right)  $.\footnote{The definition of
a bifunctor also insures the associativity of composition so that
schematically: $\hom\circ(\operatorname*{het}\circ\hom)=(\hom\circ
\operatorname*{het})\circ\hom$.}%

\begin{center}
\includegraphics[
height=2.2632in,
width=2.3869in
]%
{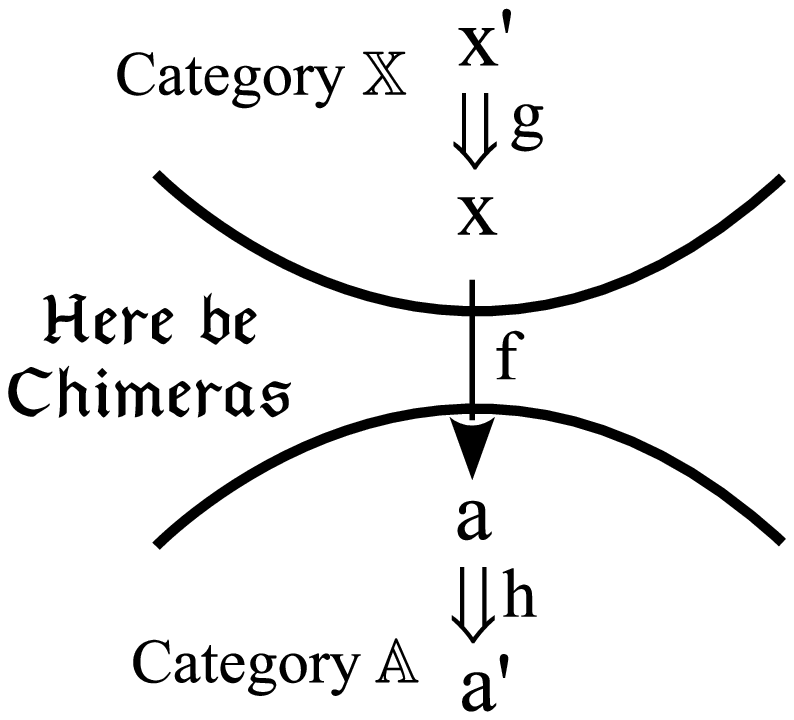}%
\end{center}

\begin{center}
Figure 15: Composition of a het with a hom on either end
\end{center}

This is all perfectly analogous to the use of $\operatorname*{Hom}$-functors
to define the composition of homs. Since both homs and hets (e.g., injection
of generators into a group) are common morphisms used in mathematical
practice, both types of bifunctors formalize standard mathematical machinery.

\subsection{Chimeras in the wilds of mathematical practice}

The homs-only orientation may go back to the original conception of category
theory "as a continuation of the Klein Erlanger Programm, in the sense that a
geometrical space with its group of transformations is generalized to a
category with its algebra of mappings." \cite[p. 237]%
{eilenberg-macl:gentheory} While chimeras do not appear in the orthodox
"ontological zoo" of category theory, they abound in the wilds of mathematical
practice. In spite of the reference to "Working Mathematician" in the title of
Mac\thinspace Lane's text \cite{maclane:cwm}, one might seriously doubt that
any working mathematician would give, say, the universal mapping property of
free groups using the "device" of the underlying set functor $U$ instead of
the traditional description given in the left representation diagram (which
does not even mention $U$) as can be seen in most any non-category-theoretic
text that treats free groups. For instance, consider the following description
in Nathan Jacobson's text \cite[p. 69]{jacobson:algebra1}.

\begin{quotation}
\noindent To summarize: given the set $X=\{x_{1},...,x_{r}\}$ we have obtained
a map $x_{i}\rightarrow\bar{x}_{i}$ of $X$ into a group $FG^{\left(  r\right)
}$ such that if $G$ is any group and $x_{i}\rightarrow a_{i}$, $1\leq i\leq r$
is any map of $X$ into $G$ then we have a unique homomorphism of $FG^{\left(
r\right)  }$ into $G$, making the following diagram commutative:
\end{quotation}

\begin{center}
$%
\begin{array}
[c]{ccc}%
X &  & \\
\downarrow & \searrow & \\
FG^{\left(  r\right)  } & \Longrightarrow & G
\end{array}
$.
\end{center}

\noindent In Jacobson's diagram, only the $FG^{\left(  r\right)
}\Longrightarrow G$ morphism is a group homomorphism; the vertical and
diagonal arrows are called "maps" and are set-to-group hets so it is the
diagram for a left representation.\footnote{We modified Jacobson's diagram
according to our het-hom convention for the arrows. Similar examples of hets
can be found in the\ Mac\thinspace Lane-Birkhoff's text \cite{maclane:algebra}%
.}

\subsection{Hets as "homs" in a collage category}

The notion of a homomorphism is so general that hets can always be recast as
"homs" in a larger category variously called a \textit{directly connected
category} \cite[p. 58]{pareigis:cats-functors} (since Pareigis calls the het
bifunctor a "connection"), a \textit{cograph} category \cite{shulman:cograph},
or, more colloquially, a \textit{collage} category (since it combines quite
different types of objects and morphisms into one category in total disregard
of any connection to the Erlangen Program). The \textit{collage category} of a
het bifunctor $\operatorname*{Het}:\mathbb{X}^{op}\times\mathbb{A}\rightarrow
Sets$, denoted $\mathbb{X}\bigstar^{\operatorname*{Het}}\mathbb{A}$ \cite[p.
96]{lurie:highertopos}, has as objects the disjoint union of the objects of
$\mathbb{X}$ and $\mathbb{A}$. The \textit{homs} of the collage category are
defined differently according to the two types of objects. For $x$ and
$x^{\prime}$ objects in $\mathbb{X}$, the homs $x\Rightarrow x^{\prime}$ are
the elements of $\operatorname*{Hom}_{\mathbb{X}}\left(  x,x^{\prime}\right)
$, the hom bifunctor for $\mathbb{X}$, and similarly for objects $a$ and
$a^{\prime}$ in $\mathbb{A}$, the homs $a\Rightarrow a^{\prime}$ are the
elements of $\operatorname*{Hom}_{\mathbb{A}}\left(  a,a^{\prime}\right)  $.
For the different types of objects such as $x$ from $\mathbb{X}$ and $a$ from
$\mathbb{A}$, the "homs" $x\Rightarrow a$ are the elements of
$\operatorname*{Het}\left(  x,a\right)  $ and there are no homs $a\Rightarrow
x$ in the other direction in the collage category.

Does the collage category construction show that "hets" are unnecessary in
category theory and that homs suffice? Since all the information given in the
het bifunctor has been repackaged in the collage category, any use of hets can
always be repackaged as a use of the "$\mathbb{X}$-to-$\mathbb{A}$ homs" in
the collage category $\mathbb{X}\bigstar^{\operatorname*{Het}}\mathbb{A}$. In
any application, like the previous example of the universal mapping property
(UMP) of the free-group functor as a left representation, one must distinguish
between the two types of objects and the three types of "homs" in the collage category.

Suppose in Jacobson's example, one wanted to "avoid" having the different
"maps" and group homomorphisms by formulating the left representation in the
collage category formed from the category of $Sets$, the category of groups
$Grps$, and the het bifunctor, $\operatorname*{Het}:Sets^{op}\times
Grps\rightarrow Sets$, for set-to-group maps. Since the UMP does not hold for
arbitrary objects and homs in the collage category, $Sets\bigstar
^{\operatorname*{Het}}Grps$, one would have to differentiate between the
"set-type objects" $X$ and the "group-type objects" $G$ as well as between the
"mixed-type homs" in $\operatorname*{Hom}\left(  X,G\right)  $ and the
"pure-type homs" in $\operatorname*{Hom}\left(  FG^{(r)},G\right)  $. Then the
left representation UMP of the free-group functor could be formulated in the
het-free collage category $Sets\bigstar^{\operatorname*{Het}}Grps$ as follows.

\begin{quotation}
\noindent For every set-type object $X$, there is a group-type object
$F\left(  X\right)  $ and a mixed-type hom $\eta_{X}:X\Rightarrow F\left(
X\right)  $ such that for any mixed-type hom $f:X\Rightarrow G$ from the
set-type object $X$ to any group-type object $G$, there is a unique pure-type
hom $f_{\ast}:F\left(  X\right)  \Rightarrow G$ such that $f=f_{\ast}\eta_{X}$.
\end{quotation}

\noindent Thus the answer to the question "Are hets really necessary?" is
"No!"--since one can always use sufficient circumlocutions with the
\textit{different} types of "homs" in a collage category. Jokes aside, the
collage category formulation is essentially only a reformulation of the left
representation UMP using clumsy circumlocutions. Working mathematicians use
phrases like "mappings" or "morphisms" to refer to hets in contrast to
homomorphisms--and "mixed-type homs" does not seem to be improved phraseology
for hets.

There is, however, a more substantive point, i.e., the general UMPs of left or
right representations show that the hets between objects of different
categories can be represented by homs \textit{within} the codomain category or
\textit{within} the domain category, respectively. If one conflates the hets
and homs in a collage category, then the point of the representation is rather
obscured (since it is then one set of "homs" in a collage category being
represented by another set of homs in the same category).

\subsection{What about the homs-only UMPs in adjunctions?}

There is another het-avoidance device afoot in the homs-only treatment of
adjunctions. For instance, the left-representation UMP of the free-group
functor can, for each $X\in Sets$, be formulated as the natural isomorphism:
$\operatorname*{Hom}_{Grps}\left(  F\left(  X\right)  ,G\right)
\cong\operatorname*{Het}\left(  X,G\right)  $. But if we fix $G$ and use the
underlying set functor $U:Grps\rightarrow Sets$, then there is trivially the
right representation: $\operatorname*{Het}\left(  X,G\right)  \cong%
\operatorname*{Hom}_{Sets}\left(  X,U\left(  G\right)  \right)  $. Putting the
two representations together, we have the heteromorphic treatment of an
adjunction first formulated by Pareigis \cite{pareigis:cats-functors}:

\begin{center}
$\operatorname*{Hom}_{Grps}\left(  F\left(  X\right)  ,G\right)
\cong\operatorname*{Het}\left(  X,G\right)  \cong\operatorname*{Hom}%
_{Sets}\left(  X,U\left(  G\right)  \right)  $.
\end{center}

\noindent If we delete the het middle term, then we have the usual homs-only
formulation of the free-group adjunction,

\begin{center}
$\operatorname*{Hom}_{Grps}\left(  F\left(  X\right)  ,G\right)
\cong\operatorname*{Hom}_{Sets}\left(  X,U\left(  G\right)  \right)  $,
\end{center}

\noindent without any mention of hets. Moreover, the het-avoidance device of
the underlying set functor $U$ allows the UMP of the free group functor to be
reformulated with sufficient circumlocutions to avoid mentioning hets.

\begin{quotation}
\noindent For each set $X$, there is a group $F\left(  X\right)  $ and a set
hom $\eta_{X}:X\Rightarrow U\left(  F\left(  X\right)  \right)  $ such that
for any set hom $f:X\Rightarrow U\left(  G\right)  $ from the set $X$ to the
underlying set $U\left(  G\right)  $ of any group $G$, there is a unique group
hom $f_{\ast}:F\left(  X\right)  \Rightarrow G$ over in the other category
such that the set hom image $U\left(  f_{\ast}\right)  $ of the group hom
$f_{\ast}$ back in the original category satisfies $f=U\left(  f_{\ast
}\right)  \eta_{X}$.\footnote{Even the "over-and-back" formulation using two
different categories could be avoided by using the further circumlocutions of
the only pure-type homs in the single collage category.}
\end{quotation}

\begin{center}%
\begin{center}
\includegraphics[
height=1.5281in,
width=2.5054in
]%
{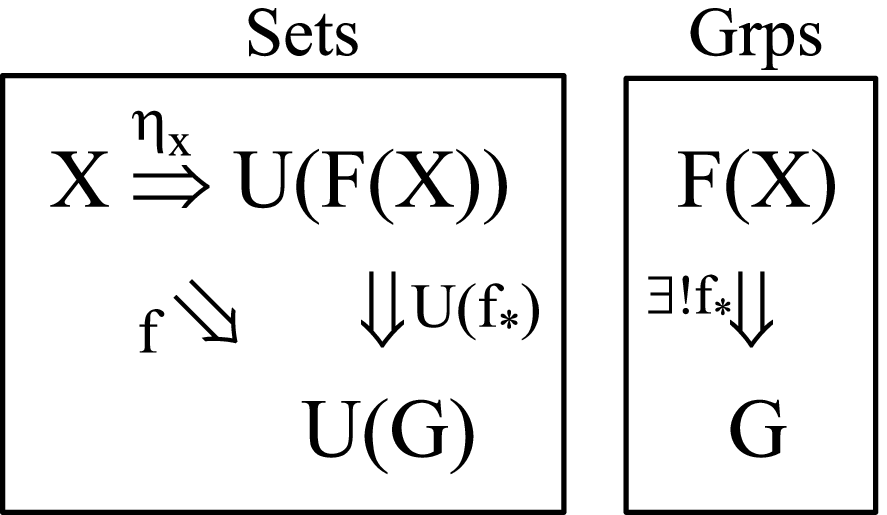}%
\end{center}

Figure 16: Over-and-back diagram for free group adjunction
\end{center}

Such het-avoidance circumlocutions have no structural significance since there
is a general adjunction representation theorem \cite[p. 147]%
{ellerman:whatisct} that \textit{all} adjoints can be represented, up to
isomorphism, as arising from the left and right representations of a het bifunctor.

\subsection{Are all UMPs part of adjunctions?}

Even though the homs-only formulation of an adjunction only ignores the
underlying hets (due to the adjunction representation theorem), is that
formulation sufficient to give all UMPs? Or are there important universal
constructions that are not either left or right adjoints?

Probably the most important example is the tensor product. The universal
mapping property of the tensor product is particularly interesting since it is
a case where the heteromorphic treatment of the UMP is forced (under one
disguise or another). The tensor product functor $\otimes:\left\langle
A,B\right\rangle \longmapsto A\otimes B$ is \textit{not} a left adjoint\ so
the usual device of using the other functor (e.g., a forgetful or diagonal
functor) to avoid mentioning hets is not available.

For $A,B,C$ modules (over some commutative ring $R$), one category is the
product category $Mod_{R}\times Mod_{R}$ where the objects are ordered pairs
$\left\langle A,B\right\rangle $ of $R$-modules and the other category is just
the category $Mod_{R}$ of $R$-modules. The values of the $\operatorname*{Het}%
$-bifunctor $\operatorname*{Het}\left(  \left\langle A,B\right\rangle
,C\right)  $ are the bilinear functions $A\times B\rightarrow C$. Then the
tensor product functor $\otimes:Mod_{R}\times Mod_{R}\rightarrow Mod_{R}$
given by $\left\langle A,B\right\rangle \longmapsto A\otimes B$ gives a left representation:

\begin{center}
$\operatorname*{Hom}_{Mod_{R}}\left(  A\otimes B,C\right)  \cong%
\operatorname*{Het}\left(  \left\langle A,B\right\rangle ,C\right)  $
\end{center}

\noindent that characterizes the tensor product. The canonical het
$\eta_{\left\langle A,B\right\rangle }:A\times B\rightarrow A\otimes B$ is the
image under the left-representation isomorphism of the identity hom
$1_{A\otimes B}$ obtained by taking $C=A\otimes B$, so we have:

\begin{center}
$%
\begin{array}
[c]{ccc}%
\left\langle A,B\right\rangle  &  & \\
^{\eta_{\left\langle A,B\right\rangle }}\downarrow^{{}} & \searrow^{f} & \\
A\otimes B & \underset{\exists!\text{ }f_{\ast}}{\Longrightarrow} & C
\end{array}
$

Left representation diagram to characterize tensor products
\end{center}

\noindent where the single arrows are the bilinear hets and the thick Arrow is
a module homomorphism within the category $Mod_{R}$.

For instance, in Mac\thinspace Lane and Birkhoff's \textit{Algebra} textbook
\cite{maclane:algebra}, they explicitly use hets (bilinear functions) starting
with the special case of an $R$-module $A$ (for a commutative ring $R$) and
then stating the universal mapping property of the tensor product $A\otimes
R\cong A$ using the left representation diagram \cite[p. 318]{maclane:algebra}%
--like any other working mathematicians. For any $R$-module $A$, there is an
$R$-module $A\otimes R$ and a canonical bilinear het $h_{0}:A\times
R\rightarrow A\otimes R$ such that given any bilinear het $h:A\times
R\rightarrow C$ to an $R$-module $C$, there is a unique $R$-module hom
$t:A\otimes R\Longrightarrow C$ such that the following diagram commutes.

\begin{center}
$%
\begin{array}
[c]{ccc}%
A\times R &  & \\
h_{0}\downarrow^{{}} & \searrow^{h} & \\
A\otimes R & \overset{\exists!t}{\Longrightarrow} & C
\end{array}
$

Left representation diagram of special case of tensor product.
\end{center}

\end{document}